\documentclass[english]{article}
\usepackage[latin1]{inputenc}
\usepackage{amsmath,amsthm,amssymb, babel}
\usepackage[notref,notcite]{}
\numberwithin{equation}{section}
\newtheorem{teo}{Theorem}

\newtheorem{lemma}{Lemma}

\def\proof{{\it Proof.}\ }
\def\endproof{\hfill $\Box$\par\vskip3mm}

\def\eq#1{(\ref{#1})}

\def\neweq#1{\begin{equation}\label{#1}}
\def\endeq{\end{equation}}

\def\phi{\varphi}

\def\RR{{\mathbb R} }

\date{}

\title{\sc  An eigenvalue problem involving a degenerate and singular elliptic operator}
\author{\sc   Mihai Mih\u ailescu$\,^{a,b}$ and Du\v san Repov\v s$\,^{c}$ \\
\small
$^a\,$Department of Mathematics, University of Craiova,  200585 Craiova,
Romania\\
\small $^b\,$Department of Mathematics, Central European
University,  1051 Budapest, Hungary\\
\small
$^c\,$Faculty of Mathematics and Physics, and Faculty of Education, University of Ljubljana,\\
\small
 POB 2964, Ljubljana, Slovenia 1001\\
\small
E-mail addresses:  {\tt mmihailes@yahoo.com}\qquad {\tt dusan.repovs@guest.arnes.si} }

\begin{document}

\maketitle \noindent{\small{\sc Abstract}. We study an eigenvalue problem involving a degenerate and singular elliptic
operator on the whole space $\RR^N$. We prove the existence of an unbounded and increasing sequence of eigenvalues. Our study generalizes
to the case of degenerate and singular operators a result of A. Szulkin and M. Willem.  
 }

\section{Introduction and main result}
The goal of this paper is to study the eigenvalue problem
\begin{equation}\label{1}
-{\rm div}(|x|^\alpha\nabla u(x))=\lambda g(x)u(x),\;\;\;\forall\;x\in\RR^N\,,
\end{equation}
where $N\geq 3$, $\alpha\in(0,2)$, $\lambda>0$ and
$g:\RR^N\rightarrow\RR$ is a function that can change sign on
$\RR^N$ satisfying the
following basic assumption\\
\smallskip

\noindent (G) $g\in L^1_{\rm loc}(\RR^N)$, $g^+=g_1+g_2\neq 0$, $g_1\in L^{\frac{N}{2-\alpha}}(\RR^N)$ and $\lim_{x\rightarrow y}|x-y|^{2-\alpha}
g_2(x)=0$, for all $y\in\RR^N$ and $\lim_{|x|\rightarrow\infty}|x|^{2-\alpha}g_2(x)=0$.
\medskip

\noindent{\bf Remark.} Note that there exists functions $h:\RR^N\rightarrow\RR$ such that $h\not\in L^{\frac{N}{2-\alpha}}(\RR^N)$ but
$h$ satisfies $\lim_{x\rightarrow y}|x-y|^{2-\alpha}h(x)=0$, for all $y\in\RR^N$ and $\lim_{|x|\rightarrow\infty}|x|^{2-\alpha}h(x)=0$.
Indeed, simple computations show that we can take $h(x)=|x|^{\alpha-2}[\log(2+|x|^{2-\alpha})]^{(\alpha-2)/N}$, if $x\neq 0$ and $h(0)=1$.
\bigskip

In the case when $\alpha=0$ problem \eq{1} becomes
\begin{equation}\label{szwi}
-\Delta u(x)=\lambda g(x)u(x),\;\;\;\forall\;x\in\RR^N\,.
\end{equation}
For this problem A. Szulkin \& M. Willem proved in \cite{SW} the
existence of an unbounded and increasing sequence of eigenvalues.
Motivated by this result on problem \eq{szwi} we consider in this
paper the natural generalization of problem \eq{szwi} given by
problem \eq{1}, obtained in the case of the presence of the
degenerate and singular potential $|x|^\alpha$ in the divergence
operator. This potential leads to a differential operator  $${\rm
div}(|x|^\alpha\nabla u(x))$$ which is degenerate and singular in
the sense that $\lim_{|x|\rightarrow 0}|x|^\alpha=0$ and
$\lim_{|x|\rightarrow\infty}|x|^\alpha=\infty$, provided that
$\alpha\in(0,2)$. Consequently, we will analyze equation \eq{1} in
the case when  the operator ${\rm div}(|x|^\alpha\nabla u(x))$ is
not strictly elliptic in the sense pointed out in D. Gilbarg \& N.
S. Trudinger \cite{GT} (see, page 31 in \cite{GT}  for the
definition of strictly elliptic operators). It follows that  some
of the techniques that can be applied in solving equations
involving strictly elliptic operators fail in this new context.
For instance some concentration phenomena may occur in the
degenerate and singular case which lead to a lack of compactness.
This is in keeping, on the one hand, with the action of the
non-compact group of dilations in $\RR^N$ and, on the other hand,
with the fact that we are looking for entire solutions for problem
\eq{1}, that means solutions defined on the whole space.
\medskip

Regarding the real-world applications of problems of type \eq{1} we remember that degenerate differential operators like the one which appears
in \eq{1} are used in the study of many physical phenomena related to equilibrium of anisotropic continuous media (see \cite{DL}). In an
appropriate context we also note that problems of type \eq{1} come also from considerations of standing waves in anisotropic Schr\"odinger
equations (see, e.g. \cite{MR}).
\medskip

A powerful tool that can be useful  when we deal with equations of type \eq{1} is the Caffarelli-Kohn-Nirenberg inequality. More exactly,
in 1984, L. Caffarelli, R. Kohn \& L. Nirenberg proved in \cite{CKN84} (see also \cite{CM1} and \cite{CM2}), in the context of some more
general inequalities, the following result:
given $p\in(1,N)$, for all $u\in C_0^\infty(\RR^N)$, there exists a positive constant $C_{a,b}$ such that
\begin{equation}\label{CKNineq}\left(\int_{\RR^N}|x|^{-bq}|u|^q\;dx\right)^{p/q}\leq C_{a,b}\int_{\RR^N}|x|^{-ap}|\nabla u|^p\;dx\,,\end{equation}
where
$$-\infty<a<\frac{N-p}{p},\;\;\; a\leq b\leq a+1,\;\;\; q=\frac{Np}{N-p(1+a-b)}\,.$$
The constant $C_{a,b}$ in inequality \eq{CKNineq} is never achieved (see the paper of F. Catrina and Z.-Q. Wang \cite{CW} for detailes).

Note that the Caffarelli-Kohn-Nirenberg inequality \eq{CKNineq} reduces to the classical Sobolev inequality (if $a = b = 0$) and to the Hardy
inequality (if $a=0$ and $b=1$). Furthermore, its utility is even more important since it implies some Sobolev and Hardy type inequalities
in the context of degenerate differential operators. More exactly, in the case when $N\geq 3$, $\alpha\in(0,2)$, $p=q=2$, $a=-\alpha/2$ and
$b=(2-\alpha)/2$ then inequality \eq{CKNineq} reads
\begin{equation}\label{hardy}
\int_{\RR^N}\frac{u^2}{|x|^{2-\alpha}}\;dx\leq C_{\frac{-\alpha}{2},\frac{2-\alpha}{2}}\int_{\RR^N}|x|^{\alpha}|\nabla u|^2\;dx,\;\;\;\forall\;
u\in C_0^\infty(\RR^N)\,.
\end{equation}
Inequality \eq{hardy} is a Hardy type inequality. The constant $C_{\frac{-\alpha}{2},\frac{2-\alpha}{2}}$ can be chosen
$\left(\frac{2}{N-2+\alpha}\right)^2$ (see, M. Willem \cite[Th\'eor\`eme 20.7]{W}). On the other hand, taking $N\geq 3$, $\alpha\in(0,2)$,
$p=2$, $q=\frac{2N}{N-2+\alpha}$, $a=-\alpha/2$, $b=0$ in \eq{CKNineq} we find that there exists a positive constant $C_\alpha:=
C_{\frac{-\alpha}{2},0}>0$ such that
the following Sobolev type inequality holds true
\begin{equation}\label{sobolev}
\left(\int_{\RR^N}|u|^{2_\alpha^\star}\;dx\right)^{2/2_\alpha^\star}\leq C_{\alpha}\int_{\RR^N}|x|^{\alpha}|\nabla u|^2\;dx,\;\;\;\forall\;
u\in C_0^\infty(\RR^N)\,,
\end{equation}
where $2_\alpha^\star=\frac{2N}{N-2+\alpha}$ plays the role of the critical Sobolev exponent in the classical Sobolev inequality.

Turning back to equation \eq{1} and taking into account the above
discussion we notice that the natural functional space where we
can analyze equation \eq{1} is the closure of $C_0^\infty(\RR^N)$
under the norm
$$\|u\|_\alpha^2=\int_{\RR^N}|x|^{\alpha}|\nabla u|^2\;dx\,.$$
Let us denote this space by $W^{1,2}_\alpha(\RR^N)$. It is easy to
see that $W^{1,2}_\alpha(\RR^N)$ is a Hilbert space with respect
to the scalar product
$$\langle u,v\rangle_\alpha=\int_{\RR^N}|x|^\alpha\nabla u\nabla v\;dx\,,$$
for all $u,\;v\in W^{1,2}_\alpha(\RR^N)$. Furthermore, according to \cite{CW} we have
$$W^{1,2}_\alpha(\RR^N)=\overline{C_0^\infty(\RR^N\setminus\{0\})}^{\|\cdot\|_\alpha}\,.$$
On the other hand, we also point out that by construction
inequalities \eq{hardy} and \eq{sobolev} hold  for all $u\in
W^{1,2}_\alpha(\RR^N)$.

We say that $\lambda>0$ is an eigenvalue of problem \eq{1} if there exists $u_\lambda\in W^{1,2}_\alpha(\RR^N)\setminus\{0\}$ such that
$$\int_{\RR^N}|x|^\alpha\nabla u_\lambda\nabla\phi\;dx=\lambda\int_{\RR^N}g(x)u_\lambda\phi\;dx\,,$$
for all $\phi\in W^{1,2}_\alpha(\RR^N)$. For each eigenvalue $\lambda>0$ we will call $u_\lambda$ in the above definition an eigenvector
corresponding to $\lambda$.

The main result of our paper is given by the following theorem:
\begin{teo}\label{T}
Assume that condition (G) is fulfilled. Then problem \eq{1} has an unbounded, increasing sequence of positive eigenvalues.
\end{teo}

\section{Proof of the main result}
The conclusion of Theorem \ref{T} will follow from the results of Lemmas \ref{l2} and \ref{l3} below.

We start by proving the following auxiliary result:
\begin{lemma}\label{l1}
Assume that condition (G) is fulfilled. Then the functional $\Lambda:W^{1,2}_\alpha(\RR^N)\rightarrow\RR$,
$$\Lambda(u)=\int_{\RR^N}g^+(x)u^2\;dx\,,$$
is weakly continuous.
\end{lemma}
\proof
$\bullet$ First, we show that $W^{1,2}_\alpha(\RR^N)\ni u\longrightarrow\int_{\RR^N}g_1(x)u^2\;dx$ is weakly continuous.

Indeed, let $\{u_n\}\subset W^{1,2}_\alpha(\RR^N)$ be a sequence converging weakly to $u\in W^{1,2}_\alpha(\RR^N)$ in $W^{1,2}_\alpha(\RR^N)$.
By \eq{sobolev} we deduce that $W^{1,2}_\alpha(\RR^N)$ is continuously embedded in $L^{2_\alpha^\star}(\RR^N)$ and consequently
$\{u_n\}$ converges weakly to $u$ in $L^{2_\alpha^\star}(\RR^N)$. It follows that $\{u_n^2\}$ converges weakly to $u^2$ in
$L^{\frac{N}{N-2+\alpha}}(\RR^N)$.

Define the operator $T:L^{\frac{N}{N-2+\alpha}}(\RR^N)\rightarrow\RR$,
$$T(\phi)=\int_{\RR^N}g_1(x)\phi\;dx\,,$$
for all $\phi\in L^{\frac{N}{N-2+\alpha}}(\RR^N)$. Undoubtedly,
$T$ is linear. Since by (G) we have $g_1\in
L^{\frac{N}{2-\alpha}}(\RR^N)$ we infer that $T$ is also
continuous. Combining that fact with the remarks considered at the
beginning   of the proof we find that
$$\lim_{n\rightarrow\infty}T(u_n)=T(u)\,,$$
in other words, $W^{1,2}_\alpha(\RR^N)\ni u\longrightarrow\int_{\RR^N}g_1(x)u^2\;dx$ is weakly continuous.


$\bullet$ Next, we verify  that $W^{1,2}_\alpha(\RR^N)\ni u\longrightarrow\int_{\RR^N}g_2(x)u^2\;dx$ is weakly continuous.
Assume again that $\{u_n\}\subset W^{1,2}_\alpha(\RR^N)$ is a sequence converging weakly to $u\in W^{1,2}_\alpha(\RR^N)$ in $W^{1,2}_\alpha(\RR^N)$
and $\epsilon>0$ is arbitrary but fixed.

By assumption (G) we deduce that there exists $R>0$ such that
$$|x|^{2-\alpha}g_2(x)\leq\epsilon,\;\;\;\forall\;x\in B_R^c(0)\,,$$
where $B_R^c(0):=\RR^N\setminus B_R(0)$ and $B_R(0)\subset\RR^N$
represents the open ball centered at the origin of radius $R$.

Since $\{u_n\}$ converges weakly to $u$ in $W^{1,2}_\alpha(\RR^N)$  we deduce that it is bounded and consequently we can define the positive constant
$$c:=\frac{2}{N-2+\alpha}\sup_n\|u_n\|_\alpha\,.$$
Inequality \eq{hardy} implies that for each $n$ we have
\begin{equation}\label{unu}
\int_{B_R^c(0)}g_2(x)u_n^2\;dx\leq\epsilon\int_{B_R^c(0)}\frac{u_n^2}{|x|^{2-\alpha}}\;dx\leq\epsilon c^2\,,
\end{equation}
and
\begin{equation}\label{doi}
\int_{B_R^c(0)}g_2(x)u^2\;dx\leq\epsilon c^2\,.
\end{equation}
Recalling again condition (G) and using a compactness argument we find that there exists a finite covering of  $\overline B_R(0)$ by closed balls
$\overline B_{r_1}(x_1)$,..., $\overline B_{r_k}(x_k)$ such that for each $j\in\{1,...,k\}$ we have
\begin{equation}\label{163276}
|x-x_j|^{2-\alpha}g_2(x)\leq\epsilon,\;\;\;\forall\;x\in\overline B_{r_j}(x_j)\,.
\end{equation}
It is easy to see that there exists $r>0$ such that for each
$j\in\{1,...,k\}$ it holds
$$|x-x_j|^{2-\alpha}g_2(x)\leq\frac{\epsilon}{k},\;\;\;\forall\;x\in\overline B_{r}(x_j)\,. $$
Defining
$$\Omega:=\cup_{i=1}^kB_{r}(x_j)$$
we have by inequality \eq{hardy} that
\begin{equation}\label{trei}
\int_{\Omega}g_2(x)u_n^2\;dx\leq\epsilon c^2\;\;\;{\rm and}\;\;\;\int_{\Omega}g_2(x)u^2\;dx\leq\epsilon c^2\,.
\end{equation}
Relation \eq{163276} implies $g_2\in L^\infty(\overline B_R(0)\setminus\Omega)$. Since $\overline B_R(0)\setminus\Omega$ is bounded we
find $g_2\in L^{\frac{N}{2-\alpha}}(\overline B_R(0)\setminus\Omega)$ and with the same arguments as in the first part of the proof we
get
\begin{equation}\label{patru}
\lim_{n\rightarrow\infty}\int_{B_R(0)\setminus\Omega}g_2(x)u_n^2\;dx=\int_{B_R(0)\setminus\Omega}g_2(x)u^2\;dx\,.
\end{equation}
Relations \eq{unu}, \eq{doi}, \eq{trei} and \eq{patru} show that
$W^{1,2}_\alpha(\RR^N)\ni u\longrightarrow\int_{\RR^N}g_2(x)u^2\;dx$ is weakly continuous.

The proof of Lemma \ref{l1} is complete. \endproof

In order to go further we consider the following minimization problem:\\
\smallskip

\noindent$(P_1)$ ${\rm minimize}_{u\in W^{1,2}_\alpha(\RR^N)}\int_{\RR^N}|x|^\alpha|\nabla u|^2\;dx$,  under restriction $\int_{\RR^N}g(x)u^2\;dx=1$.
\medskip

\begin{lemma}\label{l2}
Assume that condition (G) is fulfilled. Then problem $(P_1)$ has a solution $e_1\geq 0$. Moreover, $e_1$ is an eigenfunction of problem \eq{1}
corresponding to the eigenvalue $\lambda_1:=\int_{\RR^N}|x|^\alpha|\nabla e_1|^2\;dx$.
\end{lemma}
\proof
Consider $\{u_n\}\subset W^{1,2}_\alpha(\RR^N)$ is a minimizing sequence for $(P_1)$, i.e.
$$\int_{\RR^N}|x|^\alpha|\nabla u_n|^2\;dx\rightarrow\inf{(P_1)}\,,$$
and
$$\int_{\RR^N}g(x)u_n^2\;dx=1\,,$$
for all $n$. It follows that $\{u_n\}$ is bounded in $W^{1,2}_\alpha(\RR^N)$ and consequently there exists $u\in W^{1,2}_\alpha(\RR^N)$ such that
$\{u_n\}$ converges weakly to $u$ in $W^{1,2}_\alpha(\RR^N)$. By the weakly lower semi-continuity of the norm we deduce
$$\int_{\RR^N}|x|^\alpha|\nabla u|^2\;dx\leq\liminf_{n\rightarrow\infty}\int_{\RR^N}|x|^\alpha|\nabla u_n|^2\;dx=\inf{(P_1)}\,.$$
On the other hand, it is clear that
$$\int_{\RR^N}g^-(x)u_n^2\;dx=\int_{\RR^N}g^+(x)u_n^2\;dx-1\,,$$
for each $n$. Lemma \ref{l1} and Fatou's lemma yield
$$\int_{\RR^N}g^-(x)u^2\;dx\leq\int_{\RR^N}g^+(x)u^2\;dx-1\,,$$
or
$$1\leq \int_{\RR^N}g(x)u^2\;dx\,.$$
Define, now, $e_1=\frac{u}{(\int_{\RR^N}g(x)u^2\;dx)^{1/2}}$. It is easy to see that $\int_{\RR^N}g(x)e_1^2\;dx=1$ and
$$\int_{\RR^N}|x|^\alpha|\nabla e_1|^2\;dx=\frac{\int_{\RR^N}|x|^\alpha|\nabla u|^2\;dx}{\int_{\RR^N}g(x)u^2\;dx}\leq
\int_{\RR^N}|x|^\alpha|\nabla u|^2\;dx\leq\inf{(P_1)}\,.$$
This shows that $e_1$ is a solution of $(P_1)$. Moreover, it is easy to see that $|e_1|$ is also a solution of $(P_1)$ and consequently we can
assume that $e_1\geq 0$.

Next, for each $\phi\in W^{1,2}_\alpha(\RR^N)$ arbitrary but fixed we define $f:\RR\rightarrow\RR$ by
$$f(\epsilon)=\frac{\int_{\RR^N}|x|^\alpha|\nabla(e_1+\epsilon\phi)|^2\;dx}{\int_{\RR^N}g(x)(e_1+\epsilon\phi)^2\;dx}\,.$$
Clearly, $f$ is of class $C^1$ and $f(0)\leq f(\epsilon)$ for all $\epsilon\in\RR$. Consequently, $0$ is a minimum point of $f$ and thus,
$$f^{'}(0)=0\,,$$
or
$$\int_{\RR^N}|x|^\alpha\nabla e_1\nabla\phi\;dx\int_{\RR^N}g(x)e_1^2\;dx=\int_{\RR^N}|x|^\alpha|\nabla e_1|^2\;dx
\int_{\RR^N}g(x)e_1\phi\;dx\,.$$
Since $\phi\in W^{1,2}_\alpha(\RR^N)$ has been chosen arbitrary we deduce that the above equality holds true for each
$\phi\in W^{1,2}_\alpha(\RR^N)$. Taking into account that $\int_{\RR^N}g(x)e_1^2\;dx=1$ it follows that
$\lambda_1:=\int_{\RR^N}|x|^\alpha|\nabla e_1|^2\;dx$ is an eigenvalue of problem \eq{1} with the corresponding eigenvector $e_1$.

The proof of Lemma \ref{l2} is complete.   \endproof

In order to find other eigenvalues of problem \eq{1} we solve the minimization problems
\bigskip

\noindent$(P_n)$ ${\rm minimize}_{u\in W^{1,2}_\alpha(\RR^N)}\int_{\RR^N}|x|^\alpha|\nabla u|^2\;dx$,  under restrictions\\
$\int_{\RR^N}|x|^\alpha\nabla u\nabla e_1\;dx=...=\int_{\RR^N}|x|^\alpha\nabla u\nabla e_{n-1}\;dx=0$ and $\int_{\RR^N}g(x)u^2\;dx=1$,
\bigskip

\noindent where $e_j$ represents the solution of problem $(P_j)$, for $j\in\{1,...,n-1\}$.

\begin{lemma}\label{l3}
Assume that condition (G) is fulfilled. Then, for every $n\geq 2$ problem $(P_n)$ has a solution $e_n$. Moreover, $e_n$ is an eigenfunction of
problem \eq{1} corresponding to the eigenvalue $\lambda_n:=\int_{\RR^N}|x|^\alpha|\nabla e_n|^2\;dx$. Furthermore, $\lim_{n\rightarrow\infty}
\lambda_n=\infty$.
\end{lemma}
\proof
The existence of $e_n$ can be obtained in the same manner as in the proof of Lemma \ref{l2}, but replacing $W^{1,2}_\alpha(\RR^N)$ with the
closed linear subspace
$$X_n:=\left\{u\in W^{1,2}_\alpha(\RR^N);\;\int_{\RR^N}|x|^\alpha\nabla u\nabla e_1\;dx=...=\int_{\RR^N}|x|^\alpha\nabla u\nabla e_{n-1}\;dx=0
\;\right\}\,.$$
Then, as in Lemma \ref{l2} there exists $e_n\in X_n$ which verifies
\begin{equation}\label{cinci}
\int_{\RR^N}|x|^\alpha\nabla e_n\nabla\phi\;dx=\lambda_n\int_{\RR^N}g(x)e_n\phi\;dx,\;\;\;\forall\;\phi\in X_n\,,
\end{equation}
where $\lambda_n:=\int_{\RR^N}|x|^\alpha|\nabla e_n|^2\;dx$ and $\int_{\RR^N}g(x)e_n^2\;dx=1$.

Next, we note that for each $u\in X_n$ we have
$$\int_{\RR^N}g(x)ue_j\;dx=0,\;\;\;\forall\;j\in\{1,...,n-1\}\,,$$
and
$$\int_{\RR^N}g(x)e_je_k\;dx=\delta{j,k},\;\;\;\forall\;j,\;k\in\{1,...,n-1\}\,.$$
Consequently, for each $v\in W^{1,2}_\alpha(\RR^N)$ it holds true
$$\int_{\RR^N}g(x)\left[v-\sum_{j=1}^{n-1}\left(\int_{\RR^N}g(x)ve_j\;dx\right)e_j\right]e_k\;dx=0,\;\;\;k\in\{1,...,n-1\}\,,$$
or
$$\int_{\RR^N}|x|^\alpha\nabla\left[v-\sum_{j=1}^{n-1}\left(\int_{\RR^N}g(x)ve_j\;dx\right)e_j\right]\nabla e_k\;dx=0,\;\;\;k\in\{1,...,n-1\}\,.$$
That means
$$v-\sum_{j=1}^{n-1}\left(\int_{\RR^N}g(x)ve_j\;dx\right)e_j\in X_n\,.$$
Thus, for each $v\in W^{1,2}_\alpha(\RR^N)$ relation \eq{cinci} holds true with $\phi=v-\sum_{j=1}^{n-1}(\int_{\RR^N}g(x)ve_j\;dx)e_j$.
On the other hand,
$$0=\int_{\RR^N}|x|^\alpha\nabla e_n\nabla e_j\;dx=\lambda_j\int_{\RR^N}g(x)e_ne_j\;dx=\lambda_n\int_{\RR^N}g(x)e_ne_j\;dx\,,$$
for all $j\in\{1,...,n-1\}$. The above pieces of information yield
$$\int_{\RR^N}|x|^\alpha\nabla e_n\nabla v\;dx=\lambda_n\int_{\RR^N}g(x)e_nv\;dx,\;\;\;\forall\;v\in W^{1,2}_\alpha(\RR^N)\,,$$
i.e. $\lambda_n:=\int_{\RR^N}|x|^\alpha|\nabla e_n|^2\;dx$ is an eigenvalue of problem \eq{1} with the corresponding eigenvector $e_1$.

Next, we point out that by construction $\{e_n\}$ is an orthonormal sequence in $W^{1,2}_\alpha(\RR^N)$ and $\{\lambda_n\}$ is an increasing
sequence of positive real numbers. We show that $\lim_{n\rightarrow\infty}\lambda_n=\infty$.

Indeed, let us define the sequence $f_n:=e_n/\sqrt{\lambda_n}$. Then $\{f_n\}$ is an orthonormal sequence in $W^{1,2}_\alpha(\RR^N)$ and
$$\|f_n\|_\alpha^2=\frac{1}{\lambda_n}\int_{\RR^N}|x|^\alpha|\nabla e_n|^2\;dx=1\,,\;\;\;\forall\;n\,.$$
It means that $\{f_n\}$ is bounded in $W^{1,2}_\alpha(\RR^N)$ and consequently there exists $f\in W^{1,2}_\alpha(\RR^N)$ such that
$\{f_n\}$ converges weakly to $f$ in $W^{1,2}_\alpha(\RR^N)$.

Let $m$ an arbitrary but fixed positive integer. For each $n>m$ we have
$$\langle f_n,f_m\rangle_\alpha=0\,.$$
Passing to the limit as $n\rightarrow\infty$ we find
$$\langle f,f_m\rangle_\alpha=0\,.$$
But, the above relation holds  for each $m$ positive integer.
Consequently, we can pass to the limit as $m\rightarrow\infty$ and
we  find that
$$\|f\|_\alpha=0\,.$$
This fact implies that $f=0$ and thus, $\{f_n\}$ converges weakly
to $0$ in $W^{1,2}_\alpha(\RR^N)$. Then, by Lemma \ref{l1} we
conclude
$$\lim_{n\rightarrow\infty}\int_{\RR^N}g^+(x)f_n^2\;dx=0\,.$$
On the other hand, for each positive integer $n$ we have the estimates
$$\frac{1}{\lambda_n}=\frac{1}{\lambda_n}\int_{\RR^N}|x|^\alpha|\nabla f_n|^2\;dx=\int_{\RR^N}g(x)f_n^2\;dx\leq\int_{\RR^N}g^+(x)f_n^2\;dx\,.$$
Passing to the limit as $n\rightarrow\infty$ we find that $\lim_{n\rightarrow\infty}\lambda_n=\infty$.

The proof of Lemma \ref{l3} is complete.  \endproof
\bigskip

\noindent {\bf Acknowledgements.} This research was supported by
Slovenian Research Agency grants P1-0292-0101 and J1-2057-0101.

\end{document}